\documentclass[a4paper, 11pt]{article}

\usepackage[margin=2cm]{geometry}
\usepackage{amsmath}
\usepackage{amsfonts}
\usepackage{amssymb}
\usepackage{graphicx}
\usepackage{parskip}
\usepackage{bm}
\usepackage{todonotes}
\usepackage{caption}
\usepackage{subcaption}
\usepackage{siunitx}
\usepackage[citecolor=black, colorlinks = true, linkcolor=black]{hyperref}
\usepackage{authblk} 

\renewcommand{\r}{\bm{r}}
\newcommand{\e}{\bm{E}}
\newcommand{\gr}{\nabla}
\renewcommand{\div}{\nabla \cdot}

\newcommand{\diff}{\nabla^2}
\newcommand{\p}[1]{\partial_{#1}}
\newcommand{\ep}{\varphi}
\newcommand{\R}{\mathbb{R}}

\graphicspath{figures}

\title{A noise-robust Monte Carlo method for electric field calculations in EMC3}
\author[1,3]{William De Deyn}
\author[2]{Ruben De Wolf}
\author[3]{Vince Maes}
\author[3]{Giovanni Samaey}
\affil[1]{Institut f{\"u}r Geometrie und Praktische Mathematik, RWTH Aachen University, Germany}
\affil[2]{Department of Mechanical Engineering, KU Leuven, Belgium}
\affil[3]{Department of Computer Science, KU Leuven, Belgium}
\date{\today}

\begin{document}
	
	\maketitle

	\begin{abstract}
	EMC3 is a state-of-the-art 3D Monte Carlo code for plasma edge transport in
	stellarator configurations, but it does not yet treat the $E \times B$ drift
	self-consistently. Computing the drift requires the electric field
	$E = -\nabla \varphi$, and hence an accurate gradient of the electric
	potential $\varphi$. Because the plasma fields produced by EMC3 are
	inherently noisy, the finite difference approximation used previously
	amplifies this noise, increasingly so as the grid is refined. We extend the
	Monte Carlo Gradient Approximation method, originally developed for 1D
	Fokker--Planck equations, to a 2D setting and apply it to the electric
	potential. For an isotropic diffusion coefficient, we derive a PDE governing
	the evolution of the electric field itself, which allows $E$ to be
	approximated directly by a Monte Carlo simulation, avoiding finite
	differences altogether. A numerical experiment based on manufactured
	solutions demonstrates the accuracy of the method and shows that the
	variance grows substantially more slowly under grid refinement than for
	finite differences. The resulting formulation contains a source term
	involving the second poloidal derivative of the potential, which we neglect.
	This is admissible only when the radial scale of the plasma edge is thin
	compared to the poloidal variation scale, an ordering expected to fail near
	X-points, during detachment, and in island divertors such as W7-X. The
	present work should therefore be regarded as a proof of concept rather than
	a general electric field solver for EMC3.\\
	\\
	\noindent
	\textbf{Keywords:} Monte Carlo, gradient approximation, electric field, noise-robust
	\end{abstract}

\section{Introduction}
Nuclear fusion reactors have been a research topic for decades now. Over the years, countless machine designs were developed~\cite{Review_comm_fusion}. The most popular machine designs are the tokamak and the stellarator. The stellarator eliminates the need for an induced plasma current to generate the poloidal magnetic field~\cite{transport_and_radiation}. Instead, it relies on geometrically complex magnetic coils to induce spiraling magnetic field lines. Therefore, the mathematical models and corresponding numerical solvers have to deal with this complex 3D geometry. 

A popular code to perform simulations of the plasma edge in stellarator type devices is the EMC3-EIRENE code. The MC-MC coupled code estimates the density, momentum, and energy of the plasma. The plasma edge simulations require an efficient and accurate computation of gradients. Important gradients include the temperature gradient, pressure gradient, gradient of the electric potential, and the magnetic field gradient~\cite{W7-X_gradients,rognlien_influence_1999,plasma_gradients}. Especially the gradient of the electric potential $ \ep $ plays an important role in plasma edge simulations, as it directly determines the electric field $ \e $, through the relation $ \e = -\nabla \ep $. The electric field, in turn, is essential for modeling the $ \e \times \bm{B}$ drift, which is a key mechanism for improving the predictive capabilities of the EMC3 code~\cite{kriete2023}. The main objective, therefore, is to obtain an accurate representation of the electric field $ \e $ by reliably approximating the gradient of the electric potential $ \ep $. The electric potential itself is approximated by the EMC3 code.

In previous work~\cite{de_wolf_toward_2024}, a least squares method based on a finite difference (FD) approximation of the gradient was proposed. However, the finite difference method has the downside that it amplifies the noise introduced by the Monte Carlo approximation. We illustrate this with a simple one-dimensional example. Consider the potential, given by $\varphi(t,x)$, which is approximated with a Monte Carlo method. We quantify the variance, or noise, on the finite difference approximation as follows
\begin{equation}\label{eq:intro:var_fd}
	\text{Var}\left [ \frac{\ep_i - \ep_{i-1}}{\Delta x} \right ] = \frac{1}{\Delta x^2}\text{Var} \left[ \ep_i - \ep_{i-1} \right].
\end{equation}
In Equation \eqref{eq:intro:var_fd}, the noise of the MC approximation is given by $\text{Var} \left[ \ep_i - \ep_{i-1} \right]$. It is clear that the variance of the finite difference approximation scales quadratically with the inverse of the grid size $\Delta x$. As a consequence, the gradient approximations in EMC3 become too noisy for small grid sizes. In this paper, we propose a new method that is more robust for decreasing grid sizes than the finite difference method. The paper is based on the work in~\cite{william_mastersthesis}, which introduced and discussed the Monte Carlo Gradient Approximation (MCGA) method for general 1D Fokker-Planck equations. In this paper, we extend the method to a specific 2D case.

In Section \ref{section2}, we present the electric potential equation. The current approach to calculate the electric field uses a MC method to solve the electric potential equation and then applies the FD method. Next, Section \ref{section3} introduces the Monte Carlo Gradient Approximation method where we derive a new equation that governs the evolution of the electric field directly. We additionally highlight the challenges associated with this equation. Section \ref{section4} provides a numerical experiment to illustrate the ability and accuracy of the MC method. We also analyze the effect of grid resolution on the variance of both MC and FD approximations. Finally, Section \ref{section5} summarizes the main findings and presents the conclusions. 

\section{Electric Potential Equation}\label{section2}
In this section, we derive a reduced model for the electric potential $\ep$. 
Charge conservation requires the current density $\mathbf{J}$ to be divergence-free, i.e. 
\begin{equation}\label{eq:current_balance}
	\nabla\cdot\mathbf{J}=0.
\end{equation}

The current density is decomposed into parallel and perpendicular components, 
$\mathbf{J}=\mathbf{J_{\parallel}}+\mathbf{J_\perp}$. The parallel current $\mathbf{J_{\parallel}}$ 
follows from the electron parallel momentum balance,  
\begin{equation}
	\mathbf{J_{\parallel}}=\frac{\sigma_{\parallel}}{en_e}\left(\nabla_{\parallel}p_e - en_e\nabla_{\parallel}\ep 
	+ 0.71\,n_e\nabla_{\parallel}T_e\right),
\end{equation}
where $\sigma_{\parallel}$ is the parallel conductivity, $p_e$ and $T_e$ denote the electron 
pressure and temperature, and $n_e$ is the electron density.  

The perpendicular current $\mathbf{J_\perp}$ accounts for diffusion and non-ambipolar 
convective drifts,  
\begin{equation}
	\mathbf{J_\perp}=\mathbf{J}_{\text{drift}}-\sigma_\perp\nabla_\perp\ep,
\end{equation}
with $\sigma_\perp$ representing the anomalous perpendicular conductivity.  

Substituting these expressions into Eq.~\eqref{eq:current_balance} yields a diffusion equation 
for the electric potential,  
\begin{align}
	& \div \left( -\sigma_{\perp}\nabla_{\perp}\varphi - \sigma_{\parallel} \nabla_{\parallel}\varphi\right) = S, \label{eq:final_potential_equation} \\
	& S = - \div \bm{J}_{\text{drift}} - \div \left( \frac{\sigma_{\parallel}}{e} \left( T_{e} \nabla_{\parallel} \ln(n_{e}) + 1.71 \nabla_{\parallel} T_{e}\right) \right). \label{eq:final_source}
\end{align}

For boundary conditions at the plates, Bohm boundary conditions of the form
\begin{equation}
	\frac{e\varphi}{T_{e}} = 2.8 +0.5 \ln \left( \frac{m_{i}T_{e}}{m_{e}T}\right) - \ln \left( 1 - \frac{J_{\parallel}}{J_{sat}} \right),
\end{equation}
are employed. Other boundary conditions, such as reflective boundaries, are also possible. Equation \eqref{eq:final_potential_equation} is a stationary diffusion-type equation. However, in Sections \ref{section3} and \ref{section4}, we consider a more general time-dependent diffusion equation. For clarity and simplicity, we adopt classical Dirichlet boundary conditions; however, the approach remains valid if Bohm boundary conditions are used instead. Equation \eqref{eq:final_potential_equation} serves as the starting point for deriving a closed
evolution equation for $\e$ in the next section.

\section{Monte Carlo Gradient Approximation}\label{section3}
In this section, we introduce the Monte Carlo Gradient Approximation method for computing gradients in 2D. For both the derivation and the numerical experiments, the electric potential $ \varphi $ is modeled by a time-dependent linear diffusion equation with spatially varying isotropic diffusion coefficient $ D(\r) $ in a Cartesian coordinate system. The equation is
\begin{equation}\label{mcga:electric_pot}
	\p{t} \ep(t,\r) = \div \left(D(\r)\nabla\ep(t,\r)\right),
\end{equation}
where $ \r  = \left( x, y \right)^\top $ and $ \nabla = \left( \p{x}, \p{y} \right)^\top $. The quantity of interest is the electric field $ \e = - \nabla \ep $. The idea of the Monte Carlo Gradient Approximation method is to derive an evolution equation directly for $ \e $. Taking the gradient of Equation \eqref{mcga:electric_pot} and assuming sufficient smoothness of $ \ep $ and $ D $ to interchange the order of differentiation, we obtain
\begin{align}
	\nabla \p{t} \ep(t,\r) = \nabla\left(\div \left(D(\r)\nabla \ep(t,\r)\right)\right).
\end{align}
Next, we swap the order of $ \nabla$ and $ \p{t} $, and substitute $ \e = - \gr \ep$. This results in
\begin{equation}\label{mcga:eq_with_grad}
	\p{t} \e(t,\r) =  \nabla \left( \div \left( D(\r)\e(t,\r)\right)\right).
\end{equation} 
Finally, we reformulate Equation~\eqref{mcga:eq_with_grad} back into a standard diffusion formulation
\begin{equation}\label{mcga:final_vector}
	\p{t} \e(t,\r) = \div \left( D(\r) \nabla \e(t,\r) \right) + \underbrace{\div \left(  \nabla D(\r) \otimes \e(t,\r) \right)}_{S(t,\r)}.
\end{equation}
Equation \eqref{mcga:final_vector} is a vector diffusion equation with a source term $ S(t,\r) $.

We aim to solve Equation \eqref{mcga:final_vector} with a Monte Carlo method. However, at present, we are not aware of any Monte Carlo method that can directly approximate vector quantities. Instead, we decompose Eq.~\eqref{mcga:final_vector} into scalar equations for the components $ E_{x} $ and $ E_{y}$
\begin{equation}\label{mcga:system}
\begin{cases}
\partial_t E_x(t,\r) = \nabla \cdot (D(\r) \nabla E_x(t,\r)) + \nabla \cdot(\partial_x D(\r)\bm{E}(t,\r))\\
\partial_t E_y(t,\r) = \nabla \cdot (D(\r) \nabla E_y(t,\r)) + \nabla \cdot(\partial_y D(\r) \bm{E}(t,\r))
\end{cases}.
\end{equation}
Thus, two coupled scalar equations must be solved to approximate the electric field $ \e $.

We illustrate how these equations are coupled. Consider the source term for the $x$-component
\begin{align}
\nabla \cdot(\partial_x D(\r) \bm{E}(t,\r)) & = \partial_x \left(\partial_x D(\r) E_x(t,\r) \right) + \underline{\partial_y \left(\partial_x D(\r)E_y(t,\r) \right)}.
\end{align} 
The underlined term couples different components of the electric field through the source term, which poses difficulties. The term  $ \p{y} E_{y} $ is unknown, and cannot be incorporated into the simulation in a straightforward way. Estimating $ \p{y} E_{y}$ with finite differences would reintroduce the original problem we seek to resolve. In this work, the term $ \p{y} E_{y}$ is neglected in numerical experiments, assuming that it introduces a small bias in the simulation.

In plasma edge applications, the coordinates $x$ and $y$ are identified
with the radial and poloidal direction, respectively. The steep gradients
of the plasma edge occur in the radial direction, whereas the variation
along the poloidal direction is comparatively slow. Neglecting the
underlined term amounts to assuming that this separation of scales holds,
so that the poloidal derivative of $E_y$ is small compared to the
retained radial contributions.

The approximation can be justified when the width of the scrape-off layer
(SOL) is thin compared to the poloidal length scale. This condition is
reasonably well satisfied in the main SOL of tokamaks. However, the
ordering breaks down near X-points and during detachment, even in
tokamaks. More importantly, it is likely insufficient for devices using
an island divertor, such as W7-X, where the relevant radial scale is the
island width, which is comparable to the poloidal variation scale within
the island. In this regime, the second poloidal derivative of the
potential is not negligible and the approximation may lead to significant
errors in the computed electric field and the resulting $E \times B$
drifts. Further work is required to develop a consistent approximation strategy to include these source terms in realistic geometries.

\section{Numerical Experiments}\label{section4}
In this section, we present two experiments to evaluate the accuracy and robustness of the method. We utilize the method of manufactured solutions to construct reference solutions~\cite{roache2002}. All experiments are performed in 2D on a rectangular, uniformly discretized grid.

\subsection{Numerical Verification}
As a first test, we validate the Monte Carlo method on a simple 2D problem on a square domain. We consider the following equation for the electric potential:
\begin{equation}\label{eq:experiment1:potential}
	\p{t}\ep(t,x,y) = \div \left(Dx^2 \nabla \ep(t,x,y)\right) + S(t,x,y), \qquad D \in \R^{+}.
\end{equation}
This differs from the electric potential equation in Section \ref{section2} because Eq.~\eqref{eq:experiment1:potential} is time-dependent, the diffusion coefficient is isotropic, and the source term~\eqref{eq:final_source} is neglected. To improve readability, we omit the explicit function arguments in the following equations. We employ the method of manufactured solutions to determine the source term $S(t,x,y) $. We propose an exact solution $ \ep(x,y,t) = \exp(-(x+y))\exp(-Dt) $, from which the source term follows as
\begin{equation}
	S = \left(-D + 2Dx -2Dx^2\right)\ep
\end{equation}
The corresponding exact electric field components are
\begin{align}
	& E_{x}(t,x,y) = \exp(-(x+y))\exp(-Dt),
	& E_{y}(t,x,y) = \exp(-(x+y))\exp(-Dt).
\end{align}
Differentiating Equation \eqref{eq:experiment1:potential} with respect to $ x $ and $ y $, and substituting $ -E_x = \partial_x \ep$ and $ -E_y = \partial_y \ep$, yields the governing equations for the electric field components. For $ E_{x} $, we obtain 
\begin{equation}\label{eq:experiment1:ex}
	\p{t}E_{x} + \div \left( \begin{bmatrix}
		- 2D x \\ 0
	\end{bmatrix} E_{x}\right)  = \div \left( Dx^2 \nabla E_{x} \right) + S_{E_{x}}, 
\end{equation}  
with the source term $ S_{E_{x}} $ given by
\begin{equation}\label{eq:experiment1:source_x}
	S_{E_{x}} = \left( -2D + 4Dx \right)\ep + \left( -D + 2Dx -2Dx^2 \right)E_{x} + \underline{2Dx \p{y}E_{y}}.
\end{equation}
The underlined term is again problematic, as it introduces a coupling with the derivative of the y-component. The equation for $ E_{y} $ is simpler, since the diffusion coefficient depends only on the x-direction. We have 
\begin{equation}\label{eq:experiment1:ey}
	\p{t}E_{y}  = \div \left( Dx^2 \nabla E_{y} \right) + S_{E_{y}}, 
\end{equation}
with $ S_{E_{y}} $ equal to
\begin{equation}
	S_{E_{y}} = \left( -D + 2Dx -2Dx^2 \right)E_{y}.
\end{equation}
Dirichlet boundary conditions are imposed for $ \ep, E_{x} $, and $E_{y}$ with boundary values calculated from the exact solution. The initial condition is likewise taken from the exact solution evaluated at $ t = t_{0} $. All three equations~\eqref{eq:experiment1:potential}, \eqref{eq:experiment1:ex}, and \eqref{eq:experiment1:ey}  are simulated with identical parameter values, namely 
\begin{align}
	N = \num{500000}, \quad t_{0} = 0, \quad T^* = 1, \quad \Delta t = 0.001, \quad D = 0.1,
\end{align}
where $ N $ denotes the number of particles per simulation, $ t_{0} $ the simulation start time, $ T^* $ the end time and $ \Delta t $ the time step. The domain is $ \left[ 0,1 \right] \times \left[ 0,1 \right] $ discretized with $ M_{x} = M_{y} = 15$ grid points. In the source term \eqref{eq:experiment1:source_x}, we utilize the Monte Carlo approximation of $ \varphi $. The MC approximation of $ \ep $ is computed jointly with $ E_{x} $ and $ E_{y} $, such that three equations are solved in parallel. The source terms are incorporated using an operator-splitting method~\cite{operator_splitting}. A detailed description of the treatment of source terms can be found in~\cite{william_mastersthesis}.

Regarding the underlined term in Equation~\eqref{eq:experiment1:source_x}, two approaches are considered. Since the manufactured solution is known, $ \p{y}E_{y} $ can be evaluated exactly. When an exact solution is missing, this term is set to zero under the assumption that it has a minor contribution to the dynamics. In our test case, this assumption is well-founded, given that the solution varies only marginally in the $y$-direction. In the general setting, omitting this term is justified provided that the derivatives of the electric field are sufficiently small.

\begin{figure}[hbtp]
	\centering
	\includegraphics[width=\linewidth]{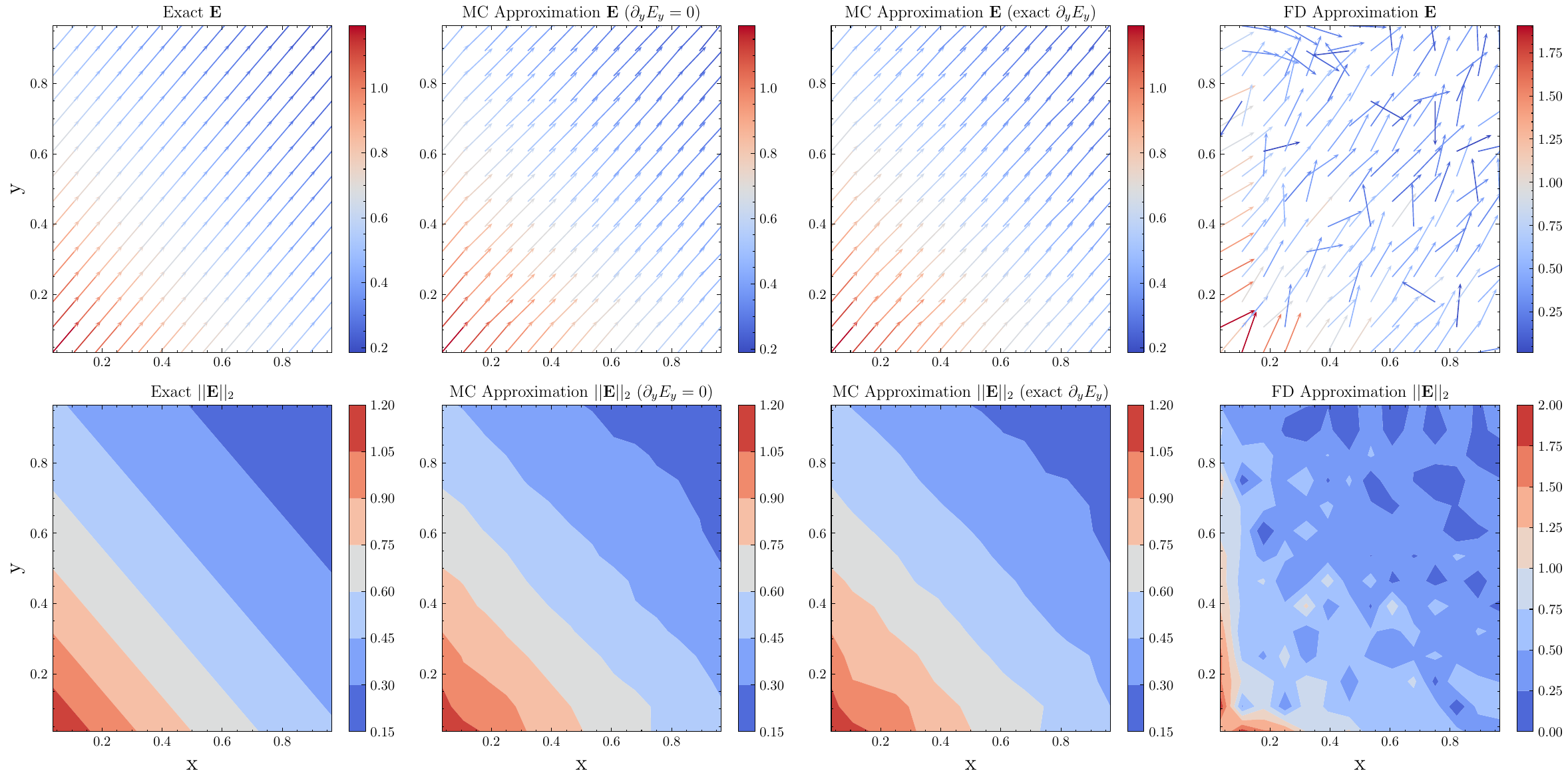}
	\captionsetup{justification=centering}
	\caption{Comparison of the MCGA and finite difference approximation of the electric field $\e$ and the norm $ \|\e \|_{2}$. For the Monte Carlo method, two cases are considered: one where the term $ \p{y} E_y $ in Eq.~\eqref{eq:experiment1:source_x} is set to zero, and another where it is given by the exact source term.}
	\label{fig:experiment1}
\end{figure}

Figure \ref{fig:experiment1} presents the electric field $ \e $ and the norm $ \|\e \|_{2}$, approximated using both the MCGA method and a first-order FD method. For the MCGA method, we additionally show the results obtained when $ \p{y}E_{y}$ is either set to zero or given by the exact source term. The Monte Carlo approximation exhibits significantly less noise compared to the finite difference approximation. Additionally, neglecting the term $ \p{y}E_{y} $ does not introduce a significant bias in the approximation, likely because in this test case the term contributes negligibly to the overall dynamics.

\begin{figure}[hbtp]
	\centering
	\includegraphics[width=\linewidth]{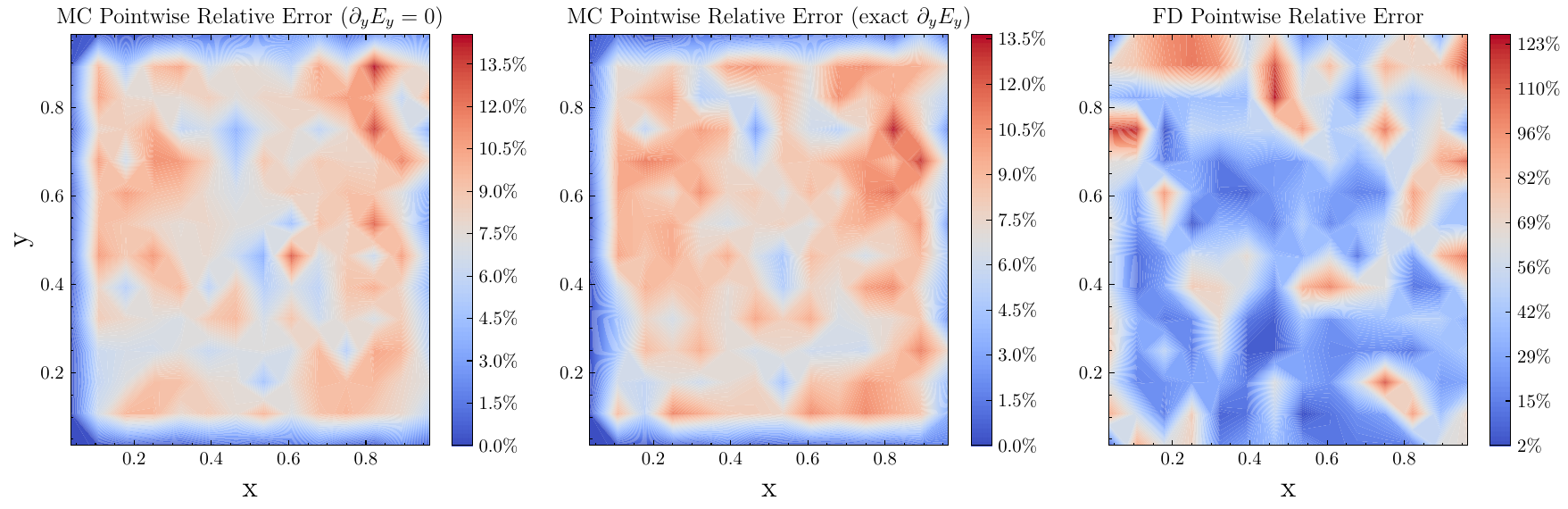}
	\captionsetup{justification=centering}
	\caption{Comparison of the relative error at each grid point for the MCGA and FD methods. As before, the MCGA method is shown both with $\partial_{y}E_{y}$ set to zero and with its exact value. Note that the colorbar scales are different between the plots.}
	\label{fig:relative-error}
\end{figure}

Figure \ref{fig:relative-error} shows the relative error of the MCGA method, both when neglecting $ \p{y} E_{y} $ and when using the exact term, as well as the relative error of the FD method. 
The FD method displays a consistently higher relative error across the domain. When $ \p{y} E_y $ is neglected, the relative error is of the same order of magnitude and the approximation remains equally accurate.

\subsection{Noise Robustness}
Next, we investigate the effect of grid resolution on the variance of the MCGA method, and compare the results with the FD method. While a Monte Carlo method is independent of the grid size, the grid resolution affects the density approximation. In 2D, for a uniform grid and equal scaling in both dimensions, the variance increases quadratically. For the FD method, we expect an increase of order 4, consistent with Eq.~\eqref{eq:intro:var_fd}. 

\begin{figure}[hbtp]
	\centering
	\includegraphics[scale = 0.35]{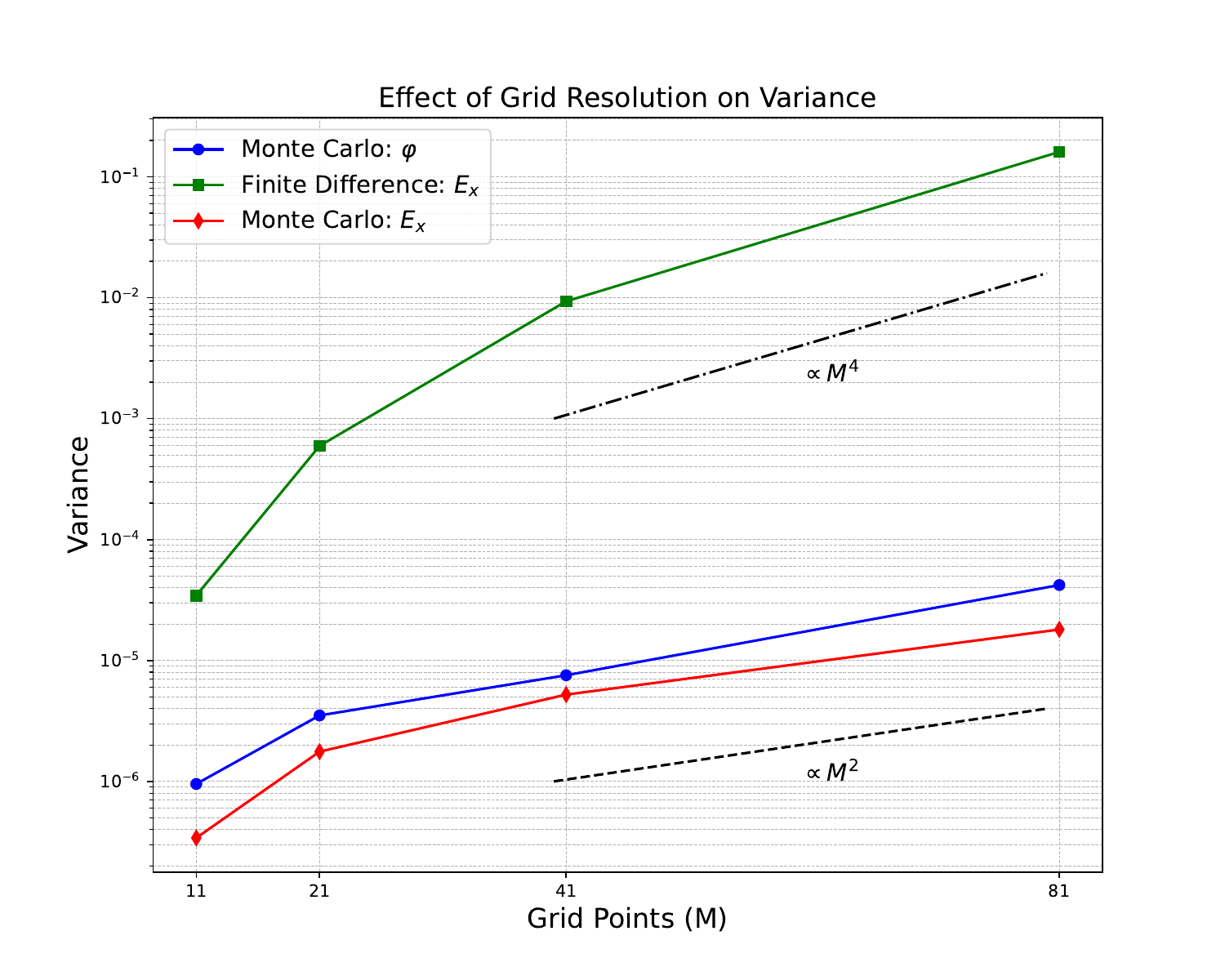}
	\captionsetup{justification=centering}
	\caption{Effect of the grid resolution on the variance. The variance of the FD approximation of $ E_{x} $ grows faster than the MC approximation of $E_{x}$. This indicates that the MC method is more robust to increases in noise as the grid size decreases.}
	\label{fig:noise_robustness}
\end{figure}

To avoid complications arising from intractable source terms, like $\partial_y E_y$, we consider the simplified diffusion problem for the electric potential
\begin{equation}\label{eq:grid:electrical_pot}
	\p{t}\ep(t,x,y) = D\diff \ep(t,x,y),
\end{equation}
with Dirichlet boundary conditions. An exact solution to Equation \eqref{eq:grid:electrical_pot} is the Gaussian 
\begin{equation}
	\ep(t,x,y) = \frac{1}{4\pi D t}\exp(-\frac{x^{2} + y^{2}}{4Dt}).
\end{equation}
The corresponding exact electric field components are
\begin{align}
	& E_{x}(t,x,y) = \frac{2x}{4\pi D t}\exp(-\frac{x^{2} + y^{2}}{4Dt}),
	& E_{y}(t,x,y) = \frac{2y}{4\pi D t}\exp(-\frac{x^{2} + y^{2}}{4Dt}),
\end{align}
which satisfy the diffusion equations
\begin{align}\label{eq:grid:electric_field}
	& \p{t} E_{x}(t,x,y) = D \diff E_{x}(t,x,y),
	& \p{t} E_{y}(t,x,y) = D \diff E_{y}(t,x,y),
\end{align}
respectively. Equations \eqref{eq:grid:electrical_pot}-\eqref{eq:grid:electric_field} are solved on the domain $ \left[ 0,2 \right] \times \left[ 0,2 \right] $ using uniform grids with different number of grid points $ M = M_{x} = M_{y} = \{11, 21, 41, 81\}$. The remaining simulation parameters are chosen as
\begin{align}
	N = \num{500000}, \quad \Delta t = 0.01, \quad t_{0} = 5, \quad T^* = 6, \quad D = 0.1.
\end{align}
The variance for each grid resolution $ M $ is computed as
\begin{equation}
	\text{Var}_{M} = \mathbb{E}\left[ \left( \hat{u}_{M} - \mathbb{E}[\hat{u}_{M}] \right)^{2}  \right],
\end{equation} 
where $ \hat{u}_{M} $ denotes the solution, calculated with either MC or FD, in the middle of the grid with $M$ grid points. We ran 20 independent simulations to approximate the expectation values with Welford's online algorithm~\cite{Welford1962NoteOA}. The results are displayed in Figure \ref{fig:noise_robustness}. We observe the theoretically predicted order of convergence. Furthermore, they indicate that the MC method exhibits greater robustness with respect to grid refinement than the FD method.

\section{Conclusion}\label{section5}
In this work, we presented the Monte Carlo Gradient Approximation method for computing the electric field. We demonstrated how to derive a new equation for the electric field from the electric potential and discussed the challenges associated with approximating the source term in this formulation. With a numerical experiment, we illustrated that the method works and achieves a lower relative error than FD, even when we set $ \p{y} E_y $ to zero. Lastly, we showed that the method achieves better robustness to noise compared to finite differences for decreasing grid sizes.

For future work, more realistic test cases need to be considered to further assess the accuracy and applicability of the Monte Carlo Gradient Approximation method. These may include cases with Bohm boundary conditions and non-rectangular grids. Such test cases are necessary to assess whether neglecting $ \p{y} E_y $ is justified. If omitting this source term introduces a significant bias, more advanced strategies for handling source terms in the electric field equations may be needed.
This is expected to be the case in island-divertor geometries, where the radial and poloidal variation scales are comparable and the thin-layer ordering underlying the present approximation no longer holds.

\section*{Author contributions}

\textbf{William De Deyn}: Conceptualization (lead),  Software implementation, Writing of the original draft. \textbf{Ruben De Wolf}: Conceptualization (support), Writing Section 2 (support), Review and editing. \textbf{Vince Maes}: Conceptualization (support), Review and editing. \textbf{Giovanni Samaey}: Conceptualization (support), Review and editing.

\section*{Acknowledgments}
The work of W.D.D. is supported by the European Union's Horizon Europe research and innovation program under the Marie Sklodowska-Curie Doctoral Network DataHyKing (Grant No. 101072546).

\bibliographystyle{abbrv}
\bibliography{references}
	
\end{document}